\documentclass[11pt]{article}
\usepackage{amsmath}
\usepackage{setspace}
\usepackage{fullpage}
\usepackage{graphicx}
\usepackage{amsthm}
\newtheorem{theorem}{Theorem}
\newtheorem{lemma}{Lemma}
\usepackage{epsf}

\begin{document}

\title{Latin squares and their defining sets}

\author{
Karola M\'esz\'aros\\
Massachusetts  Institute of Technology\\
 {\tt karola@math.mit.edu}
\\
}

\date{}

\maketitle

\begin{abstract}
\begin{small}
 A Latin square $L(n,k)$ is a square of order $n$ with its entries colored with $k$ colors so that all the entries in a row or column have different colors. Let $d(L(n,k))$ be the minimal number of colored entries of an $n \times n$ square such that there is a unique way of coloring of the yet uncolored entries in order to  obtain a Latin square $L(n, k)$. 
  In this paper we discuss the properties of $d(L(n,k))$ for $k=2n-1$ and $k=2n-2$. We give an alternate proof of the identity $d(L(n, 2n-1))=n^2-n$, which holds  for even $n$, and we establish the new result $d(L(n, 2n-2)) \geq n^2-\lfloor\frac{8n}{5}\rfloor$ and show that this bound is tight for $n$ divisible by $10$.

\end{small}
\end{abstract}

\begin{singlespace}

\section{Introduction}

The study of Latin squares has long traditions in combinatorics. The books [3, 4]  serve as reference books to the theory about Latin squares. Latin squares have various notable connections to discrete mathematics [6]. In particular, numerous problems for Latin squares naturally translate into problems in graph theory [2]. The problem of the cardinality of defining sets in graph colorings or Latin squares is one of these problems and is the focus of this paper. Results related to this problem can be found in [1, 5, 7, 8, 9].     

A \textbf{Latin square} of order $n$ is conventionally defined as a square which entries are colored with $n$ colors  so that all the entries in a row or column have different colors. In this paper we broaden the concept of a Latin square by varying the number of colors used in the coloring of the square. Let $\mathcal{L}_{n, k}=\{L \mid L$ \textit{is an} $n \times n$ \textit{square which entries are colored with} $k$ \textit{colors so that all the entries in a row or column have different colors}$\}$. 


Enumerate the columns of a $n \times n$ square from left to right and the rows from top to bottom
with 
numbers $1, 2, 3, \cdots, n$, so that we have a convenient way of referring to the entries of the square. Denote the entry in the intersection of $i^{th}$ row and $j^{th}$ column by $(i, j)$.
A coloring of a square is called a \textbf{partial coloring} if not all of the entries of the square are necessarily colored. The entries to which the partial coloring does not  assign a color are said to be uncolored.  A partial coloring \textbf{extends  to $\bf{\it{L(n, k)}}$} if there is a way to color the uncolored entries of the the given $n \times n$ square such that the  resulting entirely colored square is in $\mathcal{L}_{n, k}$. 
 A partial coloring
\textbf{uniquely extends to $L(n, k)$}  (can be uniquely extended, is uniquely extendable, etc.) if there is exactly one way to extend it to $L(n, k)$.
 A  \textbf{defining set of the $k$-coloring} for a square of order $n$ is the set of colored entries  of a partial coloring of an $n \times n$ square such that the   partial coloring uniquely extends to $L(n, k)$.
 A
defining set with minimum cardinality is called a \textbf{minimum
defining
set} and its cardinality is the \textbf{defining number}, denoted by
\(d(L(n,k))\). A defining set of the 5-coloring of a square of order 4 is presented in the figure below.

\begin {center}
\begin{tabular}{|l|l|l|l|}\hline
3&2&1&\\\hline
4&1& &2\\\hline
2&4&5&1\\\hline
 &5&2&4 \\\hline
\end{tabular}

\end{center}

For $k> 2n-1$ it is clear that  $d(L(n,k))=n^{2}$. 
However, the case \(n\le k\le 2n-1\) is not trivial. Moreover, for \(n\le k<2n-1\), 
the values of \(d(L(n,k))\) are unknown. In their paper Mahdian and Mahmoodian [1] discuss the value $d(L(n, 2n-1))$, and prove that  $d(L(n, 2n-1))=n^2-n$ for $n$ even, while $d(L(n, 2n-1))=n^2-n+1$ for  $n$ odd and greater than one.

In Section 2 of this paper we give an alternative construction showing  $d(L(n, 2n-1))=n^2-n$ for $n$ even. We exploit this construction  in Section 5 in the process of determining $d(L(n, 2n-2))$. In Section 3 we develop the tools for  proving our main result $d(L(n, 2n-2))\geq n^2-\lfloor\frac{8n}{5}\rfloor$, while Sections 4 and 5 accomplish this task as well as the verification that the bound is tight for $n$ divisible by 10.

\section{A Construction Showing  $d(L(n, 2n-1))=n^2-n$ For $n$ Even }

It has been  shown that  \(d(L(n,2n-1))=n^{2}-n\) for all even $n$ [1].
We present an alternate proof here.

It is easy to see that a  partial coloring of an $n \times n$ square $A$  cannot be extended to a $L(n, k)$ uniquely if there are more than $n$ uncolored entries of $A$. This observation immediately shows  \(d(L(n,2n-1))\geq n^{2}-n\). Therefore, if we show a partial coloring of an $n \times n$ square with exactly $n^2-n$ uncolored entries, $n$ even,  which can be uniquely extended to $L(n, 2n-1)$, then we have shown the identity  \(d(L(n,2n-1))=n^{2}-n\) for $n$ even. This is what we accomplish in the following construction. 

\textbf{Construction.}  Given an $n \times n$ square $A$, denote the union of the \(i^{th}\) row and \(j^{th}\) column by  \(u(i,j)=\{ (h, w) | (h=i\  {\rm{ and }}\   (1\le w\le n))\ {\rm{ or }}\  (w=j {\hbox{ and }}  (1\le h\le n)) \}\).
We color $A$ with $2n-1$ different colors $1, 2, \cdots, 2n-1$.

Color the entries \((i, j)\) of  \(A\),  where \(i\ne n\) and \(i>j\) with color  \(i+j \pmod{n-1}\). 
This partial coloring ensures that  the colored entries in \(u(i,i)\)
are distinct. For, the $n-2$ used colors are  \(\{i+1, i+2,\ldots, i+(i-1), (i+1)+i, (i+2)+i,\ldots, (n-1)+i\} \pmod{n-1}\), and if 
 \(i+j_1\equiv i+j_2 \pmod{n-1}\), \( (1\le j_1,j_2 \le n-1)\)  then \(j_1\equiv j_2 \pmod{n-1}\), and thus 
 \(j_1=j_2\).  Next, we color \((n, i)\) ($1\leq i\leq n-1$)
with color \(2i \pmod{n-1}\), and so  the \(n-1\) colored entries 
in \(u(i,i)\) have colors  \(\{i+1, i+2,\ldots, i+(i-1), i+i, (i+1)+i, (i+2)+i,\ldots, (n-1)+i\} \pmod{n-1} =\{1,2,3,\ldots,n-1\}\). Note that  the
\(n-1\) colored entries in the \(n^{th}\) row (\((n, i), 1\le i\le n-1)\) are distinct since 
\(2i_1\not\equiv 2i_2 \pmod{n-1}\) when \(n\) is even and \(i_1\ne i_2, 1\le i_1,i_2\le n-1\).

So far  the entries which are colored in \(u(i,i)\) contain distinct numbers (i.e. colors). 
Now we color \((j, i)\) with \(k+n\)  if \((i, j)\) has color \(k\).
This produces  a partial coloring of $A$  such that there are no numbers
(colors) on the main diagonal, but for each \(i=1,2,3,\ldots,n\) the $2n-2$ colored entries from \(u(i,i)\) are colored with the set of colors $\{1, 2, \cdots, 2n-1\} \backslash  \{n\}$.
Thus, the only possible color for the entries on the main diagonal is \(n\).

Therefore, the above  construction provides a partial coloring (with exactly $n^2-n$ colored entries)  of an $n \times n$ square, for any even $n$,  which uniquely extends  to $L(n, 2n-1)$. This completes the proof of \(d(L(n,2n-1))=n^{2}-n\) for $n$ even. 
 

\section{Properties of
Latin Squares in $\mathcal{L}_{n, 2n-2}$} 

In this section we prove that certain partial colorings of an $n \times n$ square $A$ would prohibit a unique extension of the partial coloring to $L(n, 2n-2)$. The specific features of these partial colorings enable us to prove  $d(L(n, 2n-2))\geq n^2- \lfloor\frac{8n}{5}\rfloor$  in Section 4.

To \textbf{switch} the $i^{th}$  and the $j^{th}$ row (column) of a partially colored $n \times n$ square $A$ means to switch the colors (or ``uncolors'' if some entry happens to be uncolored) of  entries $(i, k)$ and $(j, k)$ ($(k, i)$ and $(k, j)$) for all $1\leq k\leq n$.   
To $\bf{rearrange}$ $A$ means to switch some of its rows and columns. 
A partially colored square $B$ obtained from  a partially colored $n \times n$ square $A$ by some number of switching rows and columns is called  a \textbf{rearrangement} of $A$. Note that there are many possible rearrangements of $A$.
Given a partially colored $n \times n$ square $A$, and one of its rearrangements $B$, 
then  the partial coloring of $A$ uniquely extends to $L(n, k)$ if and only if the partial coloring of  $B$   uniquely extends to $L(n, k)$. 
An \textbf{available color} for an uncolored entry  \((i, j)\)  (\(1\le i,j \le n\)) in a partially colored square $L \in \mathcal{L}_{n, k}$  is one of the
\(k\) colors which
does not appear as a color of
the colored entries of $u(i,j)$. The set of available colors for uncolored entry $(i, j)$ is denoted by $a(i, j)$. Also, we denote the color of entry $(i, j)$ by $c(i, j)$.

\begin{lemma}
If the partial coloring of an $n \times n$ square  $A$  uniquely extends   to \(L(n,2n-2)\), then A has no three uncolored entries in
the same row or column.
\end{lemma}	
\begin{proof}
It suffices to prove that if the entries \((1, 1), (1, 2)\) and
\((1, 3)\) are uncolored then there is no unique extension of the partial coloring to $L(n, 2n-2)$, since the three uncolored entries can be in a row without loss of generality and we can rearrange the columns.  
We assume that the partial  coloring of \(A\)  uniquely extends
 to \(L(n,2n-2)\). 

Color all yet uncolored entries of \(A\) except
\((1, 1), (1, 2)\) and \((1, 3)\) according to the unique extension to $L(n, 2n-2)$. The partial coloring of $A$ just described (with only uncolored entries \((1, 1), (1, 2)\) and \((1, 3)\)) is also uniquely extendable to $L(n, 2n-2)$.  We prove that there exist more Latin squares in $\mathcal{L}_{n, 2n-2}$ to which we could extend the partial coloring of $A$, if there exists any such (as we supposed), obtaning the desired contradiction.

Each of the three uncolored entries \((1, 1), (1, 2)\) and \((1, 3)\) has at least two available colors, because the maximum number of different colors in the union of the row and column containing each of 
them is \(n-3+n-1=2n-4\), and we are using \(2n-2\) colors.
Suppose that some of the three uncolored entries has at least three available colors, without loss of generality let it be \((1, 1)\). Let \(e,f \in a(1, 3)\), $e\neq f$. If we set \(c(1, 3)=e\) then there is at least one remaining available color  for \((1, 2)\) and at least one of the three colors which were available at the beginning are still available for \((1, 1)\), so we can complete the coloring to a $L_1 \in \mathcal{L}_{n, 2n-2}$ this way.
Similarly, if we color  \((1, 3)\) with $f$, we can complete this partial  coloring of $A$ to a different  $L_2 \in \mathcal{L}_{n, 2n-2}$. This contradicts the unique extendability to $L(n, 2n-2)$ of $A$. 
Therefore, there are no three entries with  three or more available colors, and it remains to examine the case when   each of the three uncolored entries has  exactly two available colors.
Let  $a(1, 1)=\{a,b\}$,  $a(1, 2)=\{c,d\}$, and $a(1, 3)=\{e,f\}$, \(a\ne b, c\ne d, e\ne f\).
If some two among the sets $a(1,1)$, $a(1,2)$, $a(1,3)$ are equal, then, once given a way to color these three entries so as to obtain a $L_1 \in \mathcal{L}_{n, 2n-2}$,   we can simply permute the colors of the two entries which had the same set of available colors and obtain a coloring which gives a different 
  $L_2 \in \mathcal{L}_{n, 2n-2}$. This would contradict  unique extendability to  $L(n, 2n-2)$.
 Furthermore, if two of the sets $a(1,1)$, $a(1,2)$, $a(1,3)$  are disjoint, 
 then we can color the two entries with disjoint sets of available colors in four ways and at most  one of these ways uses both colors that were at the beginning  available for the third entry. Therefore, there would be at least three different Latin squares from  $\mathcal{L}_{n, 2n-2}$ which we could obtain. 
Thus, the  only sets of available  colors which may yield a unique extension of coloring to $L(n, 2n-2)$
are:
$a(1,1)=\{a,b\}$, $a(1,2)=\{a,c\}$, $a(1,3)=\{a,d\}$
and
$a(1,1)=\{a,b\}$, $a(1,2)=\{a,c\}$, $a(1,3)=\{b,c\}$, where  \(a,b,c,d\) are different. In the first case we have
 two different colorings of the three uncolored entries, namely set $c(1,1)=a$, $c(1,2)=c$, $c(1,3)=d$ or $c(1,1)=b$, $c(1,2)=a$, $c(1,3)=d$. In the second case color (1,1), (1, 2), (1, 3) with  \( a,c,b\) or \( b,a,c \), respectively.

Thus, if there are three uncolored entries in a row or column in a partial coloring of an $n \times n$ square, then  a unique extension to  $L(n, 2n-2)$ cannot exist.
 \end{proof}

\begin{lemma}
If the partial coloring of an $n \times n$ square  $A$ extends uniquely to  \(L(n,2n-2)\), then there are  no four uncolored entries  forming the vertices of a rectangle.
\end{lemma}	

 \begin{proof} 
 Assume the opposite and  consider  a rearrangement $B$  of $A$ with the four uncolored entries in positions $(1, 1)$, $(1, 2)$, $(2, 1)$, and $(2, 2)$. Call this set of four uncolored entries on  positions $(1, 1)$, $(1, 2)$, $(2, 1)$, and $(2, 2)$ \textit{\textbf{configuration 1}}. 
We  prove  that if there was a unique extension of a partial coloring to $\L(n, 2n-2)$, then configuration 1 could not exist in that partial coloring. This proves the lemma. 

Assume that  configuration 1 is present in $B$ and the partial coloring of $B$  uniquely extends to  $L(n, 2n-2)$. Color  the uncolored entries of $B$ except the four in configuration 1 according to this unique extension. The new partial coloring of $B$ remains uniquely extendable to $L(n, 2n-2)$.  
Note that all four  entries  $(1, 1)$, $(1, 2)$, $(2, 1)$, and $(2, 2)$
 have at least two available colors, $a, b \in a(1,1)$,  $c, d \in a(1,2)$, $e, f \in a(2,1)$, $g, h \in a(2,2)$,  \(a\ne b, c\ne d, g\ne h, e\ne f\).



Suppose  we color (1,1) with $a$ or $b$, then the  remaining three entries can be colored so as to obtain a $L \in \mathcal{L}_{n, 2n-2}$  in only one of these cases, since the partial coloring of $B$ uniquely extends to $L(n, 2n-2)$. Without loss of generality suppose that if we set $c(1,1)=a$, then  we cannot color the other three entries so as to obtain a $L \in \mathcal{L}_{n, 2n-2}$, while if we set $c(1,1)=b$, then we can uniquely extend this partial coloring to $L(n, 2n-2)$.

Next we show  \(a \in \{c,d\}\) and \(a \in \{e,f\}\).
Suppose  \(a \notin \{e,f\}\) or \(a \notin \{c,d\}\). Without loss of generality  \(a \notin \{e,f\}\).
If \(a \notin \{e,f\}\) then set $c(1,1)=a$, $c(1,2)=c$, $c(2,1)=$``whichever of \(e,f\) is not equal to the color of $(2,2)$'' and $c(2,2)=$``whichever of \(g,h\) is not equal \(c\)'' would be a possible coloring (assume \(a\ne c\)), which would contradict our previous observation  $c(1,1)\neq a$. Thus, \(a \in \{c,d\}\) and \(a \in \{e,f\}\), and without loss of generality
\(a=f\) and \(a=d\).


Furthermore,   \(\{e,c\}= \{g,h\}\). Indeed, 
suppose  \(\{e,c\}\ne \{g,h\}\). Then $c(1,1)=a$, $c(1,2)=c$, $c(2,1)=e$ and $c(2,2)=$``whichever of \(g,h\) not in \(\{e,c\}\)''
 would complete the coloring of $B$, contradicting that $c(1,1)$ cannot be $a$. 

Finally, having  $a=f=d$ and  \(\{e,c\}= \{g,h\}\) implies that we could extend $B$ to $L(n, 2n-2)$, by coloring configuration 1 in two different ways presented below.  This is the final contradiction proving the lemma.

\begin {center}
\begin{tabular}{|l|l|}\hline
\(b\) & \(a\) \\ \hline
\(a\) & \(e\) \\ \hline
\end{tabular}
\hspace{1cm}
\begin{tabular}{|l|l|}\hline
\(b\) & \(a\) \\ \hline
\(a\) & \(c\) \\ \hline
\end{tabular}
\end{center}

\end{proof}

\begin{lemma}If it is possible to uniquely extend a partial coloring to  L(n, 2n-2), then there exist no  $k, l$ such that:

$i)$ The set of available colors for $(k+1,l+1)$ is $\{a\}$, for $(k+1,l+2)$ is $\{a, b\}$, for $(k+2,l+2)$ is $\{ b, c\}$, and $(k+2,l+3)$ is also uncolored. See the left figure for illustration.  

$ii)$ The set of  available colors for $(k+3,l+2)$ is $\{a\}$, for $(k+2,l+2)$ is $\{a, b\}$,  for $(k+2,l+1)$ is $\{b, c\}$, and $(k+1,l+1)$ is also uncolored. See the  
 right figure for illustration.

\begin{center}
\begin{tabular}{|l|l|l|}\hline
\(a\) & \(a,b\)& \\ \hline
 & \(b,c\)&$\star$ \\ \hline
\end{tabular}
\hspace{1cm}
\begin{tabular}{|l|l|}\hline
$\star$& \\ \hline
 \(b,c\) & \(a,b\)\\ \hline
 & \(a\) \\ \hline
\end{tabular}
\end{center}

\end{lemma}

\begin{proof}

Since there is exactly one available color $(k+1, l+1)$ (left figrure) all entries of the $(k+1)^{th}$ row and   $(l+1)^{th}$ column must be colored.
Let \(d\) be the color of $(k+2,l+1)$. Since $a, b, c$ are available colors in the neighboring entries, \(d\ne a,b,c\). Thus, \(d\) has to
appear in
the  $(k+1)^{th}$ row or
 $(l+2)^{th}$ column since it is not an available color for $(k+1,l+2)$. However, in this case,  we cannot have
exactly one available color for
$(k+1,l+1)$ or two for $(k+2,l+2)$ (since $(k+1,l+1)$ has exactly one available color the entries of the  $(k+1)^{th}$ row and  $(l+1)^{th}$ column must be colored with different colors, and $d$ already appears in the  $(l+1)^{th}$ column, also, since $(k+2,l+2)$ 
 has exactly two available colors the entries of the  $(k+2)^{th}$ row and  $(l+2)^{th}$ column must be colored with different colors, and $d$ appears in the 
 $(k+2)^{th}$ row).
 Thus, the situation depicted in $i)$ is impossible.
The  proof of the impossibility of  $ii)$ is analogous.

\end{proof}

\begin{lemma}If  a partial coloring of an $n \times n$ square $A$  uniquely  extends to  L(n, 2n-2), then there are no five  uncolored entries in A making the configurations shown below ($\star$ placed in an entry means that the entry is uncolored, while we do not impose anything on the entries in which we did not place any mark).

\begin {center}
\begin{tabular}{|l|l|l|}\hline
$\star$&$\star$& \\\hline
 &$\star$&$\star$\\\hline
 & &$\star$\\\hline
\end{tabular}
\hspace{1cm}
\begin{tabular}{|l|l|l|}\hline
$\star$&& \\\hline
$\star$ &$\star$& \\\hline
 &$\star$ &$\star$\\\hline
 \end{tabular}

\end{center}   
\end{lemma}
 
\begin{proof}

Call the configuration of the five uncolored entries as depicted on the left of the  above figures \textit{\textbf{configuration 2}}. It suffices to prove the statement of the lemma for this of the two figures, and this is how  we  proceed.    
 Assume that the  partial coloring of $A$ contains configuration 2 and it uniquely extends to  $L(n, 2n-2)$. Suppose, without loss of generality,  that the leftmost entry of configuration 2 is on position
$(1,1)$ (it suffices to prove the statement of the lemma for any rearrangement). Color all uncolored entries of \(A\), except the entries of
configuration 2, according to the unique extension to $L(n, 2n-2)$. The obtained partial coloring  uniquely extends to $L(n, 2n-2)$.
Consider this new partial coloring. Observe that there is at least one available color for entries
$(1, 1)$  and $(3, 3)$, and there are   at least two available colors for $(1, 2)$, $(2, 2)$, and $(2, 3)$.

\textit{\textbf{First Claim.}} $(1,1)$ has exactly one available color, and  all the colored entries in the first row and  first column (that is the entries 
except $(1,1)$ and
($1,2)$) have
distinct colors.

If there were at least two available colors for
all five entries of configuration 2, then there would be at least two different ways to color these five entries and obtain a $L \in \mathcal{L}_{n, 2n-2}$, since taking any of
the (at least two) available colors for $(1,1)$  we would be able to color the other entries of configuration 2 into a $L \in \mathcal{L}_{n, 2n-2}$.
This would contradict that the partial coloring extends  uniquely  to $L(n, 2n-2)$. Thus, either (1, 1) or (3, 3) must have only a single available color. Without loss of generality, let (1, 1) have exactly one available color. This also  
implies that all the colored entries in the first row and  first column have
distinct colors.

\textit{\textbf{Second Claim.}} The cardinality of the sets  $a(1, 2)$, $a(2, 2)$, $a(2, 3)$ might  be 2 or 3, and the cardinality of $a(3, 3)$ must be 1.
 
It is easy to see that if any of the entries $(1,2)$, $(2, 2)$, $(2, 3)$, $(3,3)$  would have four available colors then the partial coloring of $A$ could not extend uniquely to $L(n, 2n-2)$. Thus the options for the cardinality of $a(1,2)$, $a(2, 2)$, $a(2, 3)$ are 2 or 3, while for  $a(3, 3)$ the options are 1, 2, or 3. 
In the following we  rule out the possibility of 2 or 3 available colors for $(3, 3)$.

If there
were three  available colors for $(3,3)$ and if there existed any way of coloring the five entries of configuration 2 so as to obtain a $L \in \mathcal{L}_{n, 2n-2}$, then    more  different $L \in \mathcal{L}_{n, 2n-2}$  could be  obtained, since only the coloring of (2,3)  may decrease  the number of available colors of (3,3), however that number can decrease only  from 3 to 2. This contradicts our assumption of unique extendability to $L(n, 2n-2)$. Furthermore, 
suppose  $a(3, 3)=\{g,h\}$, \(g\ne h\). Since the partial coloring uniquely extends to $L(n, 2n-2)$,  one of  $g$ or $h$, 
let it be \(g\) without loss of generality,  has to be  an available color  for (2,3) and
 we have to be  forced to set $c(2, 3)=g$ so as  to obtain a $L \in \mathcal{L}_{n, 2n-2}$ after the coloring of the five uncolored entries. Thus, $a(2, 3)=\{g,f\}$ and 
 $f\in a(2,2)$, such that we are forced to set $c(2, 2)=f$. This implies $a(2, 2)=\{f,e\}$ and $a(1,2)=\{e,c\}$, while 
 $a(1,1)=\{e\}$.  However, if the sets of available colors were as just described,  
then entries (1,1), (1,2),(2,2), (2,3) constitute the left figure from Lemma 3. Therefore,  there cannot be two available colors for (3, 3).  
 
Thus, there is exactly one available color for (1, 1) and (3, 3), and two or three available colors for (1,2), (2,2), (2,3). We now examine the possible cardinalities for the sets of available colors of (1,2), (2,2), (2,3) in more detail. 

\textbf{\textit{First Possibility.}}
There are   exactly two available colors for each of
(1,2), (2,2), (2,3). Since we assumed that there is only one way to color the five uncolored entries so as to obtain a $L \in \mathcal{L}_{n, 2n-2}$, we conclude that either after coloring (1,2) or (2,3)  the number of available colors for  (2,2) decreases to 1.  
If the number of available colors for (2, 2) decreased after the coloring of   
(1,2), then the number of available colors of (1,2) decreased after coloring 
(1,1) but then  entries (1, 1), (1, 2), (2, 2) and (2, 3) would constitute the left figure from Lemma 3 contradicting  the unique extendability to $L(n, 2n-2)$.   On the other hand, if the number of available colors of (2, 2) decreased after coloring (2,3), then  the number of available colors of (2, 3) decreased after coloring 
(3, 3), but then  entries (3, 3), (2, 3), (2, 2) and (1, 2) would constitute the right figure from Lemma 3 contradicting  the unique extendability to $L(n, 2n-2)$.  Thus,  it is impossible that there are exactly 1, 2, 2, 2, 1 available colors for entries (1, 1), (1, 2), (2, 2), (2, 3), (3, 3), respectively. 

\textbf{\textit{Second Possibility.}}
There are  three available colors for (1,2). It follows from our assumption of unique extendability to $L(n, 2n-2)$ that the available colors  
 must be as depicted in the figure below, colors $c_1, c_7,c_ 8$ are different, as are $c_6, c_8$ and $c_6, c_4$:
\begin {center}
\begin{tabular}{|l|l|l|}\hline
$c_1$&$c_1,c_7,c_8$& \\ \hline
 &$c_8,c_6$&$c_6,c_4$ \\ \hline
 &   & $c_4$   \\ \hline
\end{tabular}
\end{center}   

However, 
in this case,  entries  (1,2), (2,2), (2,3), (3,3) form the right figure from Lemma 3, which is a contradiction proving that the second possibility also cannot hold. If there were three available colors for (2, 3) we would have obtained a contradiction in the same manner (in that case the configuration of the left figure from Lemma 3 would appear).

There is now only one possibility remaining.

\textbf{\textit{Third Possibility.}} There are three available colors for (2,2), and two available colors for (1, 2) and (2, 3). Then, the available colors  
 must be as depicted in the figure below. The names of the colors are $c_1$, $c_2$, $c_3$, $a$, $b$. We know that $c_1 \neq c_2$, $a\neq$$c_2$, $a\neq$$c_3$, $c_2$$\neq c_3$, $b\neq$$c_3$.


\begin {center}
\begin{tabular}{|l|l|l|}\hline
$c_1$&$c_1,c_2$& \\ \hline  
 &\(a\),$c_2,c_3$&$c_3$,\(b\) \\ \hline
 &     &\(b\) \\ \hline
\end{tabular}
\end{center}   

Let $c(2,1)=k$. Since  \(c_1, c_2,c_3, a, b\) are available colors for some entries of $u(2, 1)$ we conclude that  \(k\ne c_1,c_2,c_3, a, b\). Thus,  $k \notin a(1, 2)$. Therefore, $k$ is a color of some entry in $u(1,2)$. However, it cannot appear\footnote{By ``appear'' we mean that it is a color of some entry.}
 in the $1^{st}$ row, since (1, 1) has exactly one available color, and so all the colored entries of $u(1,1)$  must be different. Thus,  \(k\) is a color of some entry in the $2^{nd}$ column. 
Analyzing $c(3,2)$  in an analogous way we obtain that $c(3,2)$ is a color of some entry in the  $2^{nd}$ row. Since there are  exactly three  available colors for
(2,2) we conclude  \(c(2,1)=c(3,2)=k\), since  otherwise there could not be \(2n-5\) different colors among the colored entries of $u(2,2)$.
Consider c(1,3)=\(l\), then \(l\ne c_1,c_2,c_3,b\) and \(l\ne k\) since there is exactly one available color for (1,1). If \(l\ne a\), then \(l\) is among the colored entries of $u(2,2)$. However, $l$ is not a color of some entry in the $2^{nd}$ row since  (2,3) has exactly two available colors,  nor  in  the $2^{nd}$ column since 
(1,2) has exactly two available colors. Therefore,  \(l=a\). Thus, \(a\ne c_1\) and \(a\ne b\).
Let \(k=c_6\), \(c_6\ne c_1,c_2,c_3,a,b\). Also,
\(a=c_4\),  \(c_4\ne\) $c_1$, $c_2$,$c_3$. Since \(l=a=c_4\), and \(l\) differs from \(b\), \(b\ne c_3,c_4,c_6\). Let \(b=c_5\) ($c_5$ could be equal $c_2$ or $c_1$).
\begin {center}
\begin{tabular}{|l|l|l|}\hline
\(c_1\) &\(c_1,c_2\)  &\(c_4\)   \\ \hline
$c_6$ &\(c_4,c_2,c_3\)&\(c_3,c_5\) \\ \hline
  &$c_6$    &\(c_5\)    \\ \hline
\end{tabular}
\end{center}

\textit{\textbf{First Observation.}} $c_2$ must appear in $u(1,1)$  since  (1,1) has exactly one available color. Since $c_2$$\in a(1, 2)$,  it cannot appear in the  $1^{st}$ row so it appears  in the $1^{st}$ column.

Denote by $\mathcal{C}_i$ the set of colors used in the colored entries of the $i^{th}$ column, and by  $\mathcal{R}_i$ the set of colors used in the colored entries of the $i^{th}$ row. 
 Let $\mathcal{A}=\mathcal{C}_1 \backslash \{c_6, c_2\}$. Note that  $\{c_6, c_2\} \subset \mathcal{C}_1$, and $c_6\neq c_2$, since $c(2,1)=c_6$ and $c_2 \in a(2,2)$.  Note that the cardinality of $\mathcal{A}$, $\mid \mathcal{A} \mid=n-3$.
Let $\mathcal{B}=\mathcal{R}_1 \backslash \{c_4\}$. 
 Note that $\mid \mathcal{B} \mid=n-3$. The intersection of any two of the sets  $\mathcal{A}$, $\mathcal{B}$, and $\{c_2,c_4,c_6\}$ is empty, since there is exactly one available color for (1, 1). 

Suppose $c_1$=$c_3$. We have $\mathcal{C}_1= \mathcal{A} \cup \{c_2, c_6\}$  and $\mathcal{R}_1=\mathcal{B} \cup \{c_4\}$.  
Analyzing the available colors for  (1,2) we get that 
 $\mathcal{C}_2=\mathcal{A} \cup \{c_6\}$.
Similarly,  considering the available colors for (2,2)
 we see that
 $\mathcal{R}_2=\{c_6\} \cup \mathcal{B}$ ($c_6\not \in \mathcal{B}$).
Since 
 $a(2,3)=\{c_1, c_5 \}$ and $a(3,3)=\{c_5\}$  we get that   $\mathcal{R}_3=\{c_1, c_6\} \cup \mathcal{B}$, and since $c_1 \not \in  \mathcal{C}_1$ column (since $c_1 \in a(1,1)$), this means that $c(3,1)\in \mathcal{B}$. However, this contradicts that all the colored entries of  $u(1,1)$ have different colors (which is necessary in order for (1, 1) to have exactly one available color). Thus, $c_1 \neq c_3$.

\textbf{\textit{Second Observation.}} $c_3$ must appear in $u(1,2)$ since   (1,2) has exactly two available colors ($c_3 \neq c_1$, $c_2$), but since $c_3 \in a(2,2)$, then $c_3$ $\not \in \mathcal{C}_2$, so $c_3 \in \mathcal{R}_1$. 

Let $\mathcal{A}$ be as defined above, and let
 $\mathcal{B}_1=\mathcal{R}_1 \backslash \{ c_4,  c_3\}$.  Note that 
$\{ c_4,  c_3\} \subset \mathcal{R}_1$  and $c_3 \neq c_4$ since  $c(1,3)=c_4$ and $c_3 \in a(2,3)$.  Also,  $\mid \mathcal{B}_1 \mid=n-4$. The intersection of any two of the sets $\mathcal{A}$, $\mathcal{B}_1$, and $\{c_2, c_3, c_4, c_6\}$ is empty since there is exactly one available color for (1, 1).
Note that $\mathcal{C}_1= \mathcal{A} \cup \{c_2, c_6\}$ and  $\mathcal{R}_1=\mathcal{B}_1 \cup \{c_3, c_4\}$.  
Analyzing the available colors for  (1,2) we get that 
$\mathcal{C}_2= \mathcal{A} \cup \{c_6\}$.
Similarly,  considering the available colors for (2,2) (and recalling $c_1\neq c_2$, $c_3$, $c_4$)
 we deduce 
 \( \mathcal{R}_1=\{c_6\} \cup \mathcal{B}_1 \cup \{c_1\}\) ($c_1$, $c_6 \not \in \mathcal{B}_1$). 
Since 
 (2,3) has exactly two colors available, once we color (3, 3) with $c_5$, the colored entries of the $3^{rd}$ column contain the set of colors $\mathcal{A} \cup \{c_4, c_5, c_2\}$. However, $c(3, 1)\in \mathcal{A} \cup \{c_2\}$, and obviously $c_5\neq c(3,1)$, since $c_5 \in a(3, 3)$, but then $c(3, 1)$ is a color in the $3^{rd}$ column. Then, however, it is impossible that there was exactly one available color for $(3, 3)$. This is the final contradiction proving Lemma 4.

 \end{proof}

\section{A Lower Bound  for  \(d(L(n, 2n-2))\)}

In this section we prove an upper bound on the number of possible uncolored entries in a partial coloring that  uniquely  extends to $L(n, 2n-2)$, giving a  lower bound for  
\(d(L(n, 2n-2))\).

\begin{theorem}If the partial coloring of an \(n\times n\) square \(A\)  extends uniquely to  
 L(n, 2n-2), then there are no more than \(\frac{8n}{5}\) uncolored entries
in
A.\end{theorem}
 
 \begin{proof} 
Call the configurations of four uncolored entries  depicted on Figures 
\ref{f9} and \ref{f10} \textit\textbf{configuration 3} and \textit{\textbf{configuration 4}}, respectively. Note that if configuration 3 or configuration 4 is present in a partial coloring of $A$ that uniquely extends to $L(n, 2n-2)$, then the rows and columns containing the uncolored entries of configuration 3 or configuration 4 cannot contain any other uncolored entries. Indeed, otherwise either there exist  three uncolored entries in a row or column, which contradicts with Lemma 1, or a rearrangement of $A$ contains configuration 2 from Lemma 4, which  contradicts the unique extendability of the partial coloring of $A$ to $L(n, 2n-2)$.

\begin{figure}
\begin {center}
\begin{tabular}{|l|l|l|}\hline
$\star$&$\star$&  \\ \hline
 &$\star$&$\star$ \\ \hline 
\end{tabular}
\end{center}
\caption{configuration 3\label{f9}}
\end{figure}
\begin{figure}
\begin {center}
\begin{tabular}{  |l|l|}\hline
  $\star$&  \\ \hline
  $\star$&$\star$ \\ \hline
   &$\star$  \\ \hline
\end{tabular}
\end{center}
\caption{configuration 4\label{f10}}
\end{figure}

To \textit{\textbf{shift a configuration to the upper left-hand-corner}} of $A$  means that if the minimal sized  rectangle into which the configuration can be fitted is of size $r_1 \times r_2$ ($r_1$ vertical, $r_2$ horizontal), and if there are no configurations already shifted to the upper left-hand-corner of $A$, then by switching rows and columns we fit the uncolored entries of the configuration   into the upper left-hand-corner rectangle with vertex positions $(1, 1)$, $(1, r_2)$, $(r_1, 1)$, $(r_1, r_2)$ (so that it stays the same configuration). If there were some  configurations already shifted to the upper left-hand-corner of $A$ with the last of them having its lowest and rightmost uncolored entry at position $(a, b)$, then fit the newly arriving configuration (with minimal sized rectangle $r_1 \times r_2$ into which it can be fit) into the rectangle with vertex positions $(a+1, b+1)$, $(a+1, b+r_2)$, $(a+r_1, b+1)$, $(a+r_1, b+r_2)$ without changing the position of rows or columns already containing shifted configurations.    

Given the partial coloring of $A$, do the following. If there is some configuration 3 or 4 present in $A$ not yet shifted to the upper left-hand corner,  shift it to the upper left-hand-corner. Repeat this until applicable. If there are no more configurations 3 or 4 present in this rearrangement of $A$ which are not yet shifted to the upper left-hand-corner, consider the uncolored entries which are in none of the configurations 3 or 4. If by switching rows and columns containing the uncolored entries which are in none of the configurations 3 or 4  some new configuration 3 or 4 can be obtained, do this. (Note that it is possible to do these rearrangings such that the configurations 3 and 4 what we already  shifted  are untouched, since there can be no other uncolored entries in the rows or columns which contain these configurations.) Shift the newly obtained configurations 3 and 4 to the upper left-hand-corner. Do this until it is possible to find some four uncolored entries which can form a new configuration 3 or 4 by switching rows and columns without changing the position of rows or columns already containing shifted configurations. Once this process is finished, let the lowest and rightmost uncolored entry which is part of a configuration 3 or 4 be $(n-n_1, n-n_2)$. The further manipulation of uncolored entries  occurs in the $n_1 \times n_2$ rectangle in the lower right hand side of the $n \times n$ square. Let this $n_1 \times n_2$ rectangle be $B$. 
Note that the only uncolored entries outside of $B$ are those of the configurations 3 and 4 since there cannot be any other uncolored entries in a row or column which contains an entry of configuration 3 or 4.

\begin{figure}
\begin {center}
\begin{tabular}{|l|l|l|l|}\hline
$\star$&$\star$& &  \\ \hline
 & &$\star$&$\star$ \\ \hline
\end{tabular}
\end{center}
\caption{configuration 5 \label{f4}}
\end{figure}
\begin{figure}
\begin {center}
\begin{tabular}{|l|l|}\hline
  $\star$&       \\ \hline
  $\star$&             \\ \hline
   &$\star$             \\ \hline
   &$\star$            \\ \hline
\end{tabular}
\end{center}
\caption{configuration 6\label{f5}}
\end{figure}

In $B$ there might be some rows  in which there are exactly two uncolored entries such that in the columns of these uncolored entries there are no more uncolored entries. Also, there might be some columns  in which there are exactly two uncolored entries such that in the rows of these uncolored entries there are no more uncolored entries. (By Lemma 1 there cannot be rows or columns  with three or more uncolored entries.) If there are more such rows or more such columns, then by switching rows and columns we can obtain \textit{\textbf{configuration 5}} or \textit{\textbf{configuration 6}} as depicted on Figures  \ref{f4} and \ref{f5}. By switching rows and columns, create as many of these configurations 5 and 6 as possible, and shift them to the upper left-hand-corner of $B$ without changing the position of rows or columns already containing shifted configurations. We stop the process when there are no more four uncolored entries not yet in a configuration 5 or 6 which could constitute either configuration 5 or 6 by means of rearrangement without changing the position of rows or columns already containing shifted configurations. 
 Once this process is finished, let the lowest and rightmost uncolored entry which is part of a configuration 5 or 6 be $(n-k, n-l)$. The further manipulation of uncolored entries occurs in the $k \times l$ rectangle in the lower right hand side of $B$ (which is the lower right hand side of the original $n \times n$ square). Call this $k \times l$ rectangle $C$. 
 Without 
loss of generality  \(l\le k\).
Note that all of the uncolored entries which are not in $C$ are entries of some configuration 3, 4, 5 or 6. 

Rectangle $C$ contains  at most one row with two uncolored entries such that the columns of these uncolored entries contain no more uncolored entries and  at most one column with two uncolored entries  such that the rows of these uncolored entries contain no more uncolored entries (since otherwise we could obtain some more of configurations 5 or 6). If there is some row with two uncolored entries such that the columns of these uncolored entries contain no more uncolored entries, switch columns so that the uncolored entries get next to each other, and let such a configuration of two such uncolored entries be named \textit{\textbf{configuration 7}}. Analogously, if there is some column with two uncolored entries such that the rows of these uncolored entries contain no more uncolored entries, switch rows so that the uncolored entries get next to each other, and let such a configuration  of two such uncolored entries be named \textit{\textbf{configuration 8}}.

\begin{figure}
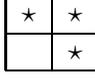

\begin {center}
\begin{tabular}{|l|l|}\hline
 $\star$&$\star$ \\ \hline
  &$\star$ \\ \hline
\end{tabular}
\end{center}
\caption{configuration 9 \label{f6}}
\end{figure} 

\begin{figure}
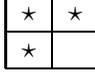

\begin {center}
\begin{tabular}{|l|l|}\hline
 $\star$&$\star$ \\ \hline
 $\star$&  \\ \hline

\end{tabular}
\end{center}
\caption{configuration 10 \label{f7}}
\end{figure}

Let the configurations of three uncolored entries from Figures \ref{f6} and \ref{f7} be \textit{\textbf{configuration 9}} and \textit{\textbf{configuration 10}}, respectively. 
By switching rows and columns obtain as many of configurations 9 and 10 as possible, without changing the position of rows or columns already containing shifted configurations. Note that there cannot be an uncolored entry which could be both in a configuration 9 and a configuration 10 depending on which switch we choose, since there can be no three uncolored entries in one row or column, and no more of configurations 3 and 4 were obtainable. From the previous and Lemma 2 it follows that there are no more uncolored entries in the rows and columns containing configurations 9 and 10. Furthermore, configurations 9 and 10 are essentially the same since switching the two columns containing them we obtain one from the other.

After introducing configurations 3, 4, 5, 6, 7, 8, 9, and 10, note that all the the uncolored entries of $A$ (more precisely of a rearrangement of $A$) must be either: an entry of configurations 3, 4, 5, 6 or in $C$. If an uncolored entry is in $C$, then it is either an entry of configurations 7, 8, 9, or 10, or it is an uncolored entry such that it is the only uncolored entry in its row and in its column. This follows from the way we defined the configurations.

Let \(l_1\) be the number of configurations 9 and 10 in $C$. 
 Let \(l_2\) be the number of columns with exactly one uncolored entry
such that the uncolored entry is the only uncolored entry in its row. 
Let
 \(L_1\) be the number of configurations 3 and 5.
Let \(L_2\) be the number of  configurations 4 and 6.

Since configurations 3, 4, 5 and 6 and the $k \times l$ rectangle $C$ must fit   inside  \(A\) without overlapping, we have the following inequalities:
\begin{equation}
3L_1+2L_2+l\le n, \end{equation}
\begin{equation}
 2L_1+3L_2+k\le n.  \end{equation}

Combining (1) and (2) we get
\begin{equation} 4(L_1+L_2)\le \frac{4}{5}(2n-k-l). \end{equation}

We distinguish four cases depending on the number of configurations 7 and 8 present in $C$. By construction this number is zero or one for each of these two configurations. 

\textit{\textbf{ Case 1.}} There is exactly one configuration 8 and no configuration 7  in $C$ ($1\leq l$, $2\leq k$). $\mathcal{U}$, the number of
 uncolored entries in  \(A\), satisfies \(\mathcal{U}= 4(L_1+L_2)+2+3l_1+l_2\).
Since
\(1+2l_1+l_2\le l\), and \(l_1\le \lfloor\frac{l-1}{2}\rfloor\) we have that 
  \(2+3l_1+l_2\le 1+l+\lfloor\frac{l-1}{2}\rfloor\).
Combining this with (3), we deduce  
\(\mathcal{U}=4(L_1+L_2)+2+3l_1+l_2\le
\frac{4}{5}(2n-k-l)+l+1+ \lfloor\frac{l-1}{2}\rfloor\), hence 
\(\mathcal{U}\le
\frac{8n}{5}- \frac{4k}{5}+\frac{l}{5}+ \lfloor\frac{l-1}{2}\rfloor+1\).

If \(l\) is even, then 
\(- \frac{4k}{5}+\frac{l}{5}+ \lfloor\frac{l-1}{2}\rfloor+1=-
\frac{4k}{5}+\frac{l}{5}+ \lfloor\frac{l}{2}\rfloor\le -\frac{4l}{5}+\frac{l}{5}+
\frac{l}{2}= \frac{-l}{10}<0\). 

If \(l\) is odd, then 
\(- \frac{4k}{5}+\frac{l}{5}+ \lfloor\frac{l-1}{2}\rfloor+1=-            
\frac{4k}{5}+\frac{l}{5}+ \lfloor\frac{l+1}{2}\rfloor\le -\frac{4l}{5}+\frac{l}{5}+
\frac{l+1}{2}= \frac{-l}{10}+\frac{1}{2}< \frac{1}{2}\). 

Thus, we obtained  $\mathcal{U}\leq \frac{8n}{5}- \frac{l}{10}<\frac{8n}{5}$ for $l$ even, while $\mathcal{U}\leq \frac{8n}{5}- \frac{l}{10}+\frac{1}{2}$  for $l$ odd. Clearly, if $l\geq 5$ and $l$ is odd, then  $\mathcal{U}\leq \frac{8n}{5}$. Therefore,  $\mathcal{U}\leq \frac{8n}{5}$ for all cases except $l\in \{1, 3\}$. Examining $n$ modulo 5, it is easily seen that  $\mathcal{U}\leq \lfloor \frac{8n}{5}\rfloor$, for all cases but   
 \(n=5m+1\) and \(l=1\), or \( n=5m+3\) and \(l\in \{1,3\}\). We  consider the remaining subcases. 

\textit{\textbf{Case 1.1.}} \(n=5m+3, l=1, k\ge 2\), or  \(n=5m+1, l=1, 
k\ge 2\). 

Since \(k>l\) then 
\begin{equation}- \frac{4k}{5}+\frac{l}{5}+ \lfloor\frac{l-1}{2}\rfloor+1=-
\frac{4k}{5}+\frac{l}{5}+ \lfloor\frac{l+1}{2}\rfloor\le -\frac{4(l+1)}{5}+\frac{l}{5}+
\frac{l+1}{2}= \frac{-l}{10}-\frac{3}{10}<0.\label{sc1}\end{equation}
 Therefore \(\mathcal{U}\le \frac{8n}{5}\).

\textit{\textbf{Case 1.2.}} \( n=5m+3, l=3,\) \( k\ge l\). 

If \(k\) 
is strictly greater than \(l\), then  (\ref{sc1}) shows \(\mathcal{U}\le \frac{8n}{5}\).
In case  \(l=k=3\),   (3) gives \(L_1+L_2\le 2m\). 
We claim that the lower right-hand corner \(3\times3\) square contains
no more then  3 uncolored entries. This holds since the rows of the entries of configuration 8 contain no other uncolored entries but the ones from configuration 8. The remaining row contains at most one uncolored entry, since we assumed there  no  configuration 7. Thus, 
\(\mathcal{U}=4(L_1+L_2)+3 \le 8m+3< 8m+4=\lfloor \frac{8(5m+3)}{5}\rfloor=\lfloor\frac{8n}{5}\rfloor\).

\textit{\textbf{Case 2.}} There is exactly one configuration 7 and no configuration 8 in $C$.
 
Since \(2+2l_1+l_2\le l\) and \(l_1\le \lfloor\frac{l-2}{2}\rfloor\), then  \( 2+3l_1+l_2= 2+2l_1+l_2+l_1\le l+\lfloor\frac{l-2}{2}\rfloor\).
Combining this with (3), we conclude
\(\mathcal{U}=4(L_1+L_2)+2+3l_1+l_2\le
\frac{4}{5}(2n-k-l)+l+ \lfloor\frac{l-2}{2}\rfloor=
\frac{8n}{5}- \frac{4k}{5}+\frac{l}{5}+ \lfloor\frac{l-2}{2}\rfloor\).

If \(l\) is even, then:
\(- \frac{4k}{5}+\frac{l}{5}+ \lfloor\frac{l-2}{2}\rfloor=-
\frac{4k}{5}+\frac{l}{5}+ \frac{l}{2}-1\le -\frac{4l}{5}+\frac{l}{5}+
\frac{l}{2}-1= \frac{-l}{10}-1< -1\).

If \( l\) is odd, then:
\(- \frac{4k}{5}+\frac{l}{5}+ \lfloor\frac{l-2}{2}\rfloor=-
\frac{4k}{5}+\frac{l}{5}+ \lfloor\frac{l-3}{2}\rfloor\le -\frac{4l}{5}+\frac{l}{5}+
\frac{l-3}{2}= \frac{-l}{10}-\frac{3}{2}<-1 \). Thus, in both cases \(\mathcal{U}<\frac{8n}{5}\).

\textbf{\textit{Case 3.}}
$C$ contains exactly one configuration 7 and one  configuration 8  ($l \geq 3$, $k \geq 3$). 
Let us put them ``up" in the lower-right-hand  rectangle of C (by switching rows and columns) as shown in the following figure:

\begin {center}
\begin{tabular}{|l|l|l|l|l|l|l|l}\hline 
 &$\star$&$\star$& & & & & \\ \hline
$\star$ && & & & & &  \\ \hline
$\star$ & &&  & & & &  \\  \hline
 & & & & & & & \\ \hline
 & & & & & & & \\ \hline
\ & & & & & & & \\ 

\end{tabular}
\end{center}

Since, 
\(3+2l_1+l_2\le l\)
and 
\(l_1\le \lfloor\frac{l-3}{2}\rfloor\), then \( 4+3l_1+l_2= 3+2l_1+l_2+l_1+1\le l+\lfloor\frac{l-3}{2}\rfloor+1\).
Combining this with (3), we obtain:
\(\mathcal{U}=4(L_1+L_2)+4+3l_1+l_2\le
\frac{4}{5}(2n-k-l)+l+ \lfloor\frac{l-3}{2}\rfloor+1=
\frac{8n}{5}- \frac{4k}{5}+\frac{l}{5}+ \lfloor\frac{l-3}{2}\rfloor+1\)

If \( l\) is even, then:
\(- \frac{4k}{5}+\frac{l}{5}+ \lfloor\frac{l-3}{2}\rfloor+1=-
\frac{4k}{5}+\frac{l}{5}+ \frac{l}{2}-1\le -\frac{4l}{5}+\frac{l}{5}+
\frac{l}{2}-1= \frac{-l}{10}-1< -1\).

If \( l\) is odd, then:
\(- \frac{4k}{5}+\frac{l}{5}+ \lfloor\frac{l-3}{2}+1\rfloor=-
\frac{4k}{5}+\frac{l}{5}+ \frac{l-1}{2}\le -\frac{4l}{5}+\frac{l}{5}+
\frac{l-1}{2}= \frac{-l}{10}-\frac{1}{2}<0 \). Thus, in both cases 
 \(\mathcal{U}<\frac{8n}{5}\).

\textit{\textbf{Case 4.}}
C contains none of configurations 7 and 8.
Then,
\(2l_1+l_2\le l\)  
and 
\(l_1\le \lfloor\frac{l}{2}\rfloor\). Hence, \( 3l_1+l_2= 2l_1+l_2+l_1\le l+\lfloor\frac{l}{2}\rfloor\).
Combining this with (3) we obtain:
\(\mathcal{U}=4(L_1+L_2)+3l_1+l_2\le
\frac{4}{5}(2n-k-l)+l+ \lfloor\frac{l}{2}\rfloor=
\frac{8n}{5}- \frac{4k}{5}+\frac{l}{5}+ \lfloor\frac{l}{2}\rfloor\)

If \( l\) is even, then:
\(- \frac{4k}{5}+\frac{l}{5}+ \lfloor\frac{l}{2}\rfloor=-
\frac{4k}{5}+\frac{l}{5}+ \frac{l}{2}\le -\frac{4l}{5}+\frac{l}{5}+
\frac{l}{2}= \frac{-l}{10}\le 0\), where the last inequality is an equality if and only if \(l=0\).

If \( l\) is odd, then:
\(- \frac{4k}{5}+\frac{l}{5}+\lfloor\frac{l}{2}\rfloor=-
\frac{4k}{5}+\frac{l}{5}+ \frac{l-1}{2}\le -\frac{4l}{5}+\frac{l}{5}+
\frac{l-1}{2}= \frac{-l}{10}-\frac{1}{2}<0 \). Thus, in both cases 
\(\mathcal{U}\leq \frac{8n}{5}\).

Therefore,  $\mathcal{U}\leq \frac{8n}{5}$ for all $n$. This proves that the number of uncolored entries in an $n \times n$ square is at most $\lfloor\frac{8n}{5}\rfloor$, if the partial coloring uniquely extends to $L(n, 2n-2)$. 

 \end{proof}

Theorem 1 immediately implies $d(L(n, 2n-2))\geq n^2-\lfloor\frac{8n}{5}\rfloor$.

\section{Example and Construction}

In this section we provide a partial coloring  for $n=5$ and $n$ divisible by 10  
with exactly  $n^2-\lfloor\frac{8n}{5}\rfloor$ colored entries, such that the presented partial coloring  uniquely extends  to $L(n, 2n-2)$. This combined with Theorem 1 proves  
$d(L(n, 2n-2))= n^2-\lfloor\frac{8n}{5}\rfloor$ for $n=5$ and $n$ divisible by 10.

\noindent {\textbf{The Construction for \(n=5\).}}

\noindent A partial coloring uniquely extending to $L(5, 8)$ of the $5 \times 5$ square with 8 uncolored entries (denoted by $\star$ in the figure) using 8 different colors $1, 2, 3, \cdots, 8$ is:

\begin {center}
\begin{tabular}{|l|l|l|l|l|}\hline
$\star$ &$\star$&7 &8 &4\\\hline
3 &$\star$&$\star$ &1 &8\\\hline
2 &6 & 5 &7 &$\star$\\\hline
5 &7 &6&$\star$ &$\star$ \\\hline
6 &5 &2&$\star$ & 3\\\hline
\end{tabular}
\end{center}

The available colors for the uncolored entries are 
\begin{center}
\begin{tabular}{|l|l|l|l|l|}\hline
\bf{1} &\bf{1,2,3}&7 &8 &4 \\ \hline
3 &\bf{2,4}&\bf{4} &1 &8 \\ \hline
2 &6 & 5 &7 &\(\bf{1}\) \\ \hline
5 &7 &6&\(\bf{2,3,4}\) &\(\bf{1,2}\)  \\ \hline
6 &5 &2&\(\bf{4}\) & 3 \\ \hline
\end{tabular}
\end{center}

We fill in the entries which have exactly one available color:
\begin{center}
\begin{tabular}{|l|l|l|l|l|}\hline
\it{1} &\bf{2,3}&7 &8 &4 \\ \hline
3 &\bf{2}&\it{4} &1 &8 \\ \hline
2 &6 & 5 &7 &\(\it{1}\) \\ \hline
5 &7 &6&\(\bf{2,3}\) &\(\bf{2}\)  \\ \hline
6 &5 &2&\(\it{4}\) & 3 \\ \hline
\end{tabular}
\end{center} 

We repeat the process: 
\begin{center}
\begin{tabular}{|l|l|l|l|l|}\hline
\it{1} &\bf{3}&7 &8 &4 \\ \hline
3 &\it{2}&\it{4} &1 &8 \\ \hline  
2 &6 & 5 &7 &\(\it{1}\) \\ \hline 
5 &7 &6&\(\bf{3}\) &\(\it{2}\)  \\ \hline
6 &5 &2&\(\it{4}\) & 3 \\ \hline
\end{tabular}
\end{center}

We uniquely get:

\begin {center}
\begin{tabular}{|l|l|l|l|l|l|}\hline
1 &3&7 &8 &4 \\ \hline
3 &2&4&1 &8 \\ \hline
2 &6 & 5 &7 &1 \\ \hline
5 &7 &6&3 &2  \\ \hline
6 &5 &2&4 & 3 \\ \hline 
\end{tabular}
\end{center} 

The construction for $n$, where $n$ is divisible by 10,  uses the construction for the $5 \times 5$ square, as well as our construction from Section 2. 

\noindent \textbf{The construction for $n$ divisible by 10.}

\noindent Color an \(\frac{n}{5}\times\frac{n}{5}\) square, $\frac{n}{5}$ even, with  \(\frac{2n}{5}-1\) colors as described in Section 2. Color the entries on the main diagonal as well, with color $\frac{n}{5}$. Define a correspondence $f$ from  $[\frac{2n}{5}-1]$ to the set of subsets of $[2n-2]$  as follows:
 $f(\frac{n}{5})=\{1, 2, 3, 4, 5, 6, 7, 8\}$, 
 $f(i)=\{8+5(i-1)+1, 8+5(i-1)+2, 8+5(i-1)+3, 8+5(i-1)+4, 8+5(i-1)+5\}$, for $1\leq i< \frac{n}{5}$, and 
 $f(i)=\{8+5(i-2)+1, 8+5(i-2)+2, 8+5(i-2)+3, 8+5(i-2)+4, 8+5(i-2)+5\}$, for $\frac{n}{5} <i \leq \frac{2n}{5}-1$.
Note that ${\bigcup}^{\frac{2n}{5}-1}_{i=1} f(i)=[2n-2]$ and $f(i) \bigcap f(j)=\emptyset$, for $i\ne j$, $i, j \in \{1, 2, 3, \cdots, \frac{2n}{5}-1\}$.

Now we  show how to color the $n \times n$ square, $n$ divisible by 10, with $2n-2$ different colors $1, 2, 3, \cdots, 2n-2$, so that there are exactly $\frac{8n}{5}$ uncolored entries and there is exactly one way of coloring the uncolored entries so as to obtain a $L \in \mathcal{L}_{n, 2n-2}$. 

 Given the $n \times n$ square,  divide it into \(5\times 5\) squares so that we  have exactly
\(\frac{n}{5}\times\frac{n}{5}\) little squares. See the figure below for illustration in the case \(n=10\):

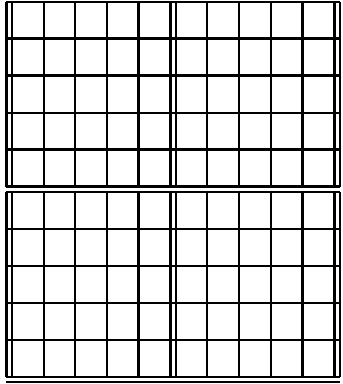
\begin {center}
\begin{tabular}{||l|l|l|l|l||l|l|l|l|l||}\hline\hline
& & & & & & & & & \\ \hline
& & & & & & & & & \\ \hline
& & & & & & & & & \\ \hline
& & & & & & & & & \\ \hline
& & & & & & & & & \\ \hline\hline
& & & & & & & & & \\ \hline
& & & & & & & & & \\ \hline
& & & & & & & & & \\ \hline
& & & & & & & & & \\ \hline
& & & & & & & & & \\ \hline\hline
\end{tabular}
\end{center}

Think of the correspondence defined above as a way of going from a coloring of the \(\frac{n}{5}\times\frac{n}{5}\) square to the coloring of the $n \times n$ square. Namely, if square $(i, j)$ in the  \(\frac{n}{5}\times\frac{n}{5}\) square was colored with color $k\neq \frac{n}{5}$ and $f(k)=\{k_1, k_2, k_3, k_4, k_5\}$, then the $5 \times 5$ square, with vertex positions $(5(i-1)+1, 5(j-1)+1)$,  $(5(i-1)+1, 5(j-1)+5)$,  $(5(i-1)+5, 5(j-1)+1)$,  $(5(i-1)+5, 5(j-1)+5)$ in the $n \times n$ square is colored ``cyclically'', as shown in the figure below:

\begin {center}
\begin{tabular}{|l|l|l|l|l|l|}\hline
$k_1$ &$k_2$&$k_3$ &$k_4$ &$k_5$ \\ \hline
$k_2$ &$k_3$&$k_4$&$k_5$ &$k_1$ \\ \hline
$k_3$ &$k_4$ & $k_5$ &$k_1$ &$k_2$ \\ \hline
$k_4$  &$k_5$  &$k_1$ &$k_2$  &$k_3$   \\ \hline
$k_5$  &$k_1$  &$k_2$ &$k_3$  & $k_4$  \\ \hline 
\end{tabular}
\end{center} 

Also, the \(5\times5 \) squares on the main diagonal  are colored with eight colors,
\(\{1,2,3,4,5,6,7,8\}\) as in the partial coloring of  the \(5\times 5 \) square in the construction for $n=5$; if some entry was uncolored in the construction, let it be uncolored here, too. 

Note that the partial coloring of the $n \times n$ square just described uniquely  extends  to $L(n, 2n-2)$. Also, there are exactly $\frac{8n}{5}$ uncolored entries in this partial coloring. Since we have proven in Section 4 that  $d(L(n, 2n-2))\geq n^2 -\lfloor\frac{8n}{5}\rfloor$,  the previous construction proves that $d(L(n, 2n-2))=n^2 -\frac{8n}{5}$ for all $n$, $n$ divisible by 10.

We illustrate the above construction for $n=10$:

Color the $\frac{n}{5} \times \frac{n}{5}$ square with   $\frac{2n}{5}-1$ 
colors:
\begin {center}
\begin{tabular}{|l|l|}\hline
2& 1  \\ \hline
3 &2 \\ \hline
\end{tabular}
\end{center}

Correspondence $f$ is given by:
\(f(2)=\{1,2,3,4,5,6,7,8\}\),
\(f(1)=\{9,10,11,12,13\}\),
\(f(3)=\{14,15,16,17,18\}\).

The partial coloring of $10 \times 10$ square is:

\begin {center}
\begin{tabular}{||l|l|l|l|l||l|l|l|l|l||}\hline\hline
$\star$ &$\star$&7 &8 &4&9&10&11&12&13 \\ \hline
3 &$\star$&$\star$ &1 &8 &10&11&12&13&9\\ \hline
2 &6 & 5 &7 &$\star$&11&12&13&9&10 \\ \hline
5 &7 &6&$\star$ &$\star$ &12&13&9&10&11 \\ \hline
6 &5 &2&$\star$ & 3 &13&9&10&11&12\\ \hline\hline
14&15&16&17&18&$\star$ &$\star$&7 &8 &4 \\ \hline
15&16&17&18&14&3 &$\star$&$\star$ &1 &8 \\ \hline
16&17&18&14&15&2 &6 & 5 &7 &$\star$ \\ \hline
17&18&14&15&16&5 &7 &6&$\star$ &$\star$  \\ \hline
18&14&15&16&17&6 &5 &2&$\star$ & 3 \\ \hline \hline
\end{tabular}
\end{center}

Analogously to the  case \(n=5\), we get the unique extension to $L(n, 2n-2)$:

\begin {center}
\begin{tabular}{||l|l|l|l|l||l|l|l|l|l||}\hline\hline
1 &3&7 &8 &4&9&10&11&12&13 \\ \hline 
3 &2&4 &1 &8 &10&11&12&13&9\\ \hline
2 &6 & 5 &7 &1&11&12&13&9&10 \\ \hline
5 &7 &6&3 &2 &12&13&9&10&11 \\ \hline
6 &5 &2&4 & 3 &13&9&10&11&12\\ \hline\hline
14&15&16&17&18&1 &3&7 &8 &4 \\ \hline 
15&16&17&18&14&3 &2&4 &1 &8 \\ \hline
16&17&18&14&15&2 &6 & 5 &7 &1 \\ \hline
17&18&14&15&16&5 &7 &6&3 &2  \\ \hline
18&14&15&16&17&6 &5 &2&4 & 3 \\ \hline \hline
\end{tabular}
\end{center}

\section{Summary}
We determined a lower  bound for the defining number, \(d(L(n,2n-2)) \ge n^{2}-\lfloor\frac{8n}{5}\rfloor\), and  exhibited a construction for \(n=5\) and \(n\) divisible by 10  when this bound
is reached.  We anticipate that a suitable lower bound for $d(L(n,2n-3))$ can be obtained  by recognizing prohibited configurations in a manner similar to that  in our paper.

\section*{Acknowledgments}
I am grateful to  Roya Beheshti Zavareh for her support and advice as my mentor  at the Research Science Institute 2000, MIT,  where this research was done. 
I thank Rado\v s  Radoi\v ci\'c for his comments after reading a draft of this paper.   
I am thankful to  Professor Hartley Rogers, 
and  The Center for Excellence in Education, for
helping and supporting the research.
 

\begin{singlespace}

\end{singlespace}

\end{singlespace}
\end{document}